\documentclass[12pt]{article}
\def\Dj{\hbox{D\kern-.73em\raise.30ex\hbox{-}
\raise-.30ex\hbox{}}}
\def\dj{\hbox{d\kern-.33em\raise.80ex\hbox{-}
\raise-.80ex\hbox{\kern-.40em}}}
\usepackage{epsfig}
\usepackage{amsmath,amsthm,amsfonts,amssymb,amscd}
\usepackage[center]{titlesec}
\allowdisplaybreaks[4]
\newtheorem {Lemma}{Lemma}
\newtheorem {Theorem} {Theorem}

\newtheorem {Corollary}{Corollary}
\newenvironment {Proof} {\noindent {\bf Proof.}}{\quad $\blacksquare$\par\vspace{3mm}}

\begin{document}
\baselineskip=0.30in

\vspace*{40mm}

\begin{center}

{\Large \bf ON DISTANCE SPECTRAL RADIUS AND DISTANCE ENERGY OF GRAPHS}

\vspace{10mm}

{\large \bf Bo Zhou\:\!$^a$ \ and \  Aleksandar Ili\' c\:\!$^b$}

\vspace{6mm}

\baselineskip=0.20in

$^a${\it Department of Mathematics, South China Normal University,\\
Guangzhou 510631, P. R. China}\\
{\rm e-mail:} {\tt zhoubo@scnu.edu.cn} \\[2mm]
$^b${\it Faculty of Sciences and Mathematics, University of
Ni\v s, \\
Vi\v segradska 33, 18000 Ni\v s, Serbia} \\
{\rm e-mail:} {\tt aleksandari@gmail.com}

\vspace{6mm}

(Received November ?, 2009)

\end{center}

\vspace{6mm}

\baselineskip=0.27in

\noindent {\bf Abstract }

\vspace{4mm}

For a connected graph, the distance spectral radius is the largest
eigenvalue of its distance matrix, and the distance energy is
defined as the sum of the absolute values of the eigenvalues of its
distance matrix. We establish lower and upper  bounds for the
distance spectral radius of graphs and bipartite graphs, lower
bounds for the distance energy of graphs, and characterize the
extremal graphs. We also discuss upper bounds for the distance
energy.

\baselineskip=0.30in


\section*{\normalsize 1. INTRODUCTION}

Let $G$ be a connected graph with vertex set $\{v_1, v_2,\dots,
v_n\}$\,. The distance between vertices $v_i$ and $v_j$ of $G$\,,
denoted by $d_{ij}$\,, is defined to be the length (i.~e., the number
of edges) of the shortest path from $v_i$ to $v_j$\,. The distance
matrix of $G$\,, denoted by ${\bf D}(G)$\,, is the $n\times n$ matrix
whose $(i,j)$-entry is equal to $d_{ij}$\,, $i,j=1,2,\ldots,n$\,, (see
\cite{CDS,FH}). Note that $d_{ii}=0$\,, $i=1,2,\ldots, n$\,. The
eigenvalues of ${\bf D}(G)$ are said to be the $D$-eigenvalues of
$G$\,. Since ${\bf D}(G)$ is a real symmetric matrix, the
$D$-eigenvalues are real and can be
ordered in non-increasing order, $\rho _{1}\geq \rho _{2}\geq \cdots \geq \rho _{n}$\,.
The distance spectral radius of $G$ is the largest $D$-eigenvalue $\rho_1$ and denoted by $\rho(G)$\,.

Balaban et al. \cite{BaCi91} proposed the use of $\rho (G)$ as a
molecular descriptor, while in \cite{GuMe98} it was successfully
used to infer the extent of branching and model boiling points of
alkanes.  In \cite{Zh07}, the author gave lower and upper bounds for
$\rho (G)$ when $G$ is a tree. In \cite{ZhTr07}, the authors
provided lower and upper bounds for $\rho (G)$ when $G$ is a
connected graph in terms of the number of vertices, the sum of the
squares of the distances between all unordered pairs of vertices and
the sum of the distances between a given vertex and all other
vertices, and the Nordhaus--Gaddum--type result for $\rho (G)$\,,
see also \cite{ZhTr07-2} for more results on  $\rho (G)$\,. A survey
on the properties of $\rho (G)$ may be found in \cite{ZT-review}.
Recently, Das \cite{Da09} obtained lower and upper bounds for the
distance spectral radius of a connected bipartite graph and
characterize those graphs for which these bounds are best possible.
Note that earlier study of the eigenvalues of the distance matrix
may be found in [10--14]. 

The distance energy of a connected graph $G$ is defined in
\cite{IGV} as
$$
DE(G)=\sum_{i=1}^{n}\left\vert \rho _{i}\right\vert\,.
$$
Lower and upper bounds for distance energy have been obtained in [15--19]. 
For more recent results on
$DE$ see [20--23]. 

In this paper, we establish lower and upper  bounds for the distance
spectral radius of graphs and bipartite graphs, lower bounds for the
distance energy of graphs, and characterize the extremal graphs. We
also discuss upper bounds for the distance energy.


\section*{\normalsize 2.  PRELIMINARIES}

Let $K_n$ be the complete graph with $n$ vertices. Let $P_n$ be the
path with $n$ vertices. Let $K_{p,q}$ be the complete bipartite
graph with $p$ vertices in one partite set and $q$ vertices in the
other partite set. Let $d_i$ be the degree of vertex $v_i$ of the
graph $G$\,. A graph is semi--regular if it is bipartite and all
vertices in the same partite set have the same degree.

Let $G$ be a connected graph.  Let $v_rv_s$ be an edge of $G$ such
that $G - v_rv_s$ is also connected. Then $d_{ij}(G)\leq
d_{ij}(G-v_rv_s)$ for all $i,j\in V(G)$\,. Moreover, $1 = d_{rs}(G)
< d_{rs}(G-v_rv_s)$ and thus by the Perron--Frobenius theorem, we
conclude that
\begin{equation}
\label{eq-deledge} \rho (G) < \rho (G - v_rv_s)\,.
\end{equation}
Similarly, for two nonadjacent vertices $v_r$ and $v_s$\,,
\begin{equation}
\label{eq-addedge} \rho (G+v_rv_s) < \rho (G)\,.
\end{equation}

Let $G$ be a graph with $n$ vertices.
Let ${\bf A}(G)$ be the adjacency matrix of the graph $G$\,.
The eigenvalues $\lambda_i$\,, $i=1, 2, \dots, n$\,, of $G$
are the eigenvalues of its adjacency matrix ${\bf A}(G)$\,, and they can be ordered as
$\lambda _{1}\geq \lambda _{2}\geq \cdots \geq \lambda _{n}$\,
 \cite{CDS}. If $G$ is $r$-regular, then $\lambda _1=r$\,. Let $\overline{G}$ be the complement of $G$\,.
Denote by
${\bf J}_n$ the all $1$'s $n\times n$ matrix and by ${\bf I}_n$ the
identity matrix of order $n$\,.  If the diameter of $G$ is at most two,
then ${\bf D}(G)=2\ {\bf J}_n-2\ {\bf I}_n-{\bf
A}(G)={\bf J}_n-{\bf I}_n+{\bf
A}(\overline{G})$ (see \cite{CDGT})\,.

For an $r$-regular graph $G$\,, the
eigenvectors of ${\bf A}(G)$ associated to any eigenvalue not equal
to $r$ are orthogonal to the all $1$'s vector.  If $G$ is an
$r$-regular graph of diameter two, then  ${\bf D}(G)=2\ {\bf J}_n-2\ {\bf I}_n-{\bf
A}(G)$\,, and thus the
$D$-eigenvalues of $G$ are $2n-r-2$\,, $-\lambda_n-2$\,, \dots,
$-\lambda_2-2$\,,  arranged in a
non-increasing manner.

\section*{\normalsize 3.  BOUNDS FOR $\rho$ OF GENERAL GRAPHS}

In this section, we present lower and upper bound for  $\rho(G)$ of a connected graph $G$\,.

\begin{Theorem} \label{low1}
Let $G$ be a connected graph with $n$ vertices, maximum degree $\Delta_1$ and second maximum degree $\Delta_2$\,.
Then
$$
\rho (G) \geq \sqrt{(2n-2-\Delta_1)(2n-2-\Delta_2)}
$$
with equality if and only if $G$ is a regular graph with diameter less than or equal to $2$\,.
\end{Theorem}

\begin{Proof}
Let $x = (x_1, x_2, \ldots, x_n)^T$
be a Perron eigenvector of $D(G)$ corresponding to the largest eigenvalue $\rho (G)$\,, such that
$$
x_i = \min_{k \in V(G)} x_k \quad \mbox{and} \quad x_j = \min_{k \in V(G)\atop k\ne i} x_k\,.
$$
From the eigenvalue equation $\rho (G) \cdot x = D (G) \cdot x$\,,
written for the component $x_i$ we have
\begin{eqnarray*}
\rho(G) x_i &=& \sum_{k = 1}^n d_{ik} \cdot x_k\\
&\geq & d_i x_j + (n-1-d_i) 2x_j = (2n-2-d_i) x_j\,.
\end{eqnarray*}
Analogously for the component $x_j$ we have
\begin{eqnarray*}
\rho(G) x_j &=& \sum_{k = 1}^n d_{jk} \cdot x_k\\[2mm]
&\geq & d_j x_i + (n-1-d_j) 2x_i = (2n-2-d_j) x_i\,.
\end{eqnarray*}
Combining these two inequalities, it follows
$$
\rho (G) \geq \sqrt{(2n-2-d_i)(2n-2-d_j)}\geq \sqrt{(2n-2-\Delta_1)(2n-2-\Delta_2)}\,.
$$

The equality holds if and only if the diameter of $G$ is less than or equal to $2$\,,
and all coordinates $x_i$ are equal. For $d = 1$\,, we get a complete graph $K_n$\,.
For $d = 2$\,, we get
$\rho (G) x_i = d_i x_i + 2 (n - 1 - d_i) x_i$\,, and then  $\rho (G) = 2n - 2 - d_i$\,,
which means that $G$ is a regular graph. Conversely, it is easily seen that
$\rho (G) = 2n - 2 - \Delta_1$ if $G$ is a regular graph with diameter less than or equal to $2$\,.
\end{Proof}

\begin{Theorem}
Let $G$ be a connected graph with $n$ vertices, minimum degree~$\delta_1$ and second minimum degree $\delta_2$\,.
Let $d$ be the diameter of $G$\,. Then
$$
\rho (G) \leq \sqrt{\left [ dn - \frac{d (d - 1)}{2} - 1 - \delta_1(d-1) \right ] \left [ dn - \frac{d (d - 1)}{2} - 1 - \delta_2(d-1) \right ]}
$$
with equality if and only if $G$ is a regular graph with diameter less than or equal to $2$\,.
\end{Theorem}

\begin{Proof}
Let $x = (x_1, x_2, \ldots, x_n)^T$
be a Perron eigenvector of $D(G)$ corresponding to the largest eigenvalue $\rho (G)$\,, such that
$$
x_i = \max_{k \in V(G)} x_k \quad \mbox{and} \quad x_j = \max_{k \in V(G)\atop k\ne i} x_k\,.
$$
From the eigenvalue equation $\rho (G) \cdot x = D (G) \cdot x$\,,
written for the component $x_i$  we have
\begin{eqnarray*}
\rho(G) x_i &=& \sum_{k = 1}^n d_{ik} \cdot x_k \\[2mm]
&\leq & d_i x_j + 2 x_j + 3 x_j + \cdots + (d-1) x_j \\[2mm]
&& + \ d \left[n-1-d_i - (d - 2)\right] x_j\\[2mm]
&=& \left [ dn - \frac{d (d - 1)}{2} - 1 - d_i(d-1) \right ] x_j\,.
\end{eqnarray*}
Analogously for the component $x_j$ we have
\begin{eqnarray*}
\rho(G) x_j &=& \sum_{k = 1}^n d_{jk} \cdot x_k \\[2mm]
&\leq & d_j x_i + 2 x_i + 3 x_i + \cdots + (d-1) x_i\\[2mm]
&& + \ d \left[n-1-d_j - (d - 2)\right] x_i\\[2mm]
&=& \left [ dn - \frac{d (d - 1)}{2} - 1 - d_j(d-1) \right ] x_i\,.
\end{eqnarray*}
Combining these two inequalities, it follows that
\begin{eqnarray*}
\rho (G) &\leq & \sqrt{\left [ dn - \frac{d (d - 1)}{2} - 1 - d_i(d-1) \right ]
\left [ dn - \frac{d (d - 1)}{2} - 1 - d_j(d-1) \right
]}\\[2mm]
&\leq &  \sqrt{\left [ dn - \frac{d (d - 1)}{2} - 1 - \delta_1(d-1) \right ]
\left [ dn - \frac{d (d - 1)}{2} - 1 - \delta_2(d-1) \right ]}\,.
\end{eqnarray*}

The equality holds if and only if all coordinates of Perron's eigenvector are equal,
and hence, $D(G)$ has equal row sums. If the diameter of $G$ is greater than or equal to~$3$\,,
that means that for every vertex $i$\,, there is exactly one vertex $j$ that is
of distance two from $i$\,, and then  the diameter of $G$ must be smaller than $4$\,.
If the diameter of $G$ is $3$ and equality holds, then for a center vertex $s$ (with the eccentricity two),
from $\rho (G) \cdot x = D (G) \cdot x$\,,
written for the component $x_s$\,, we have
$$
\rho (G)x_s =  d_s x_s + (n-1-d_s) 2 x_s=\left [ 3n - \frac{3 (3 - 1)}{2} - 1 - d_s(3-1) \right ] x_s
$$
and then $d_s=n-2$\,, which implies that $G\cong P_4$\,. But the coordinates of Perron's eigenvector of $D(P_4)$ can not be all equal.
Therefore, in the case of equality we have that $G$ is a regular graph with diameter $d \leq 2$\,.
\end{Proof}

\section*{\normalsize 4.  BOUNDS FOR $\rho$ OF BIPARTITE GRAPHS}

The inequality (\ref{eq-deledge}) shows that the maximum distance
spectral radius will be attained for trees. Subhi and Powers in
\cite{SuPo90} proved that for $n \ge 3$ the path $P_n$ has the
maximum distance spectral radius among trees with $n$ vertices.
Stevanovi\'c and Ili\'c in \cite{StIl09} generalized this result,
and proved that among trees with $n$ vertices and fixed maximum
degree $\Delta$, the broom graph $B_{n, \Delta}$ (formed by
attaching $\Delta-1$  pendent vertices to an end vertex of the path
$P_{n-\Delta+1}$) has the maximum $\rho$-value.

The inequality (\ref{eq-addedge}) tells us  that the complete
bipartite graph $K_{p, q}$ has the minimum distance spectral radius
among connected bipartite graphs  with $p$ vertices in one partite
set and $q$ vertices in the other partite set. Following
\cite{CDGT}, the distance spectrum of the complete bipartite graph
$K_{p,q}$ consists of simple eigenvalues $p + q - 2 \pm \sqrt{p^2 -
pq + q^2}$\,, and an eigenvalue $-2$ with multiplicity $p + q -
2$\,. Let $G$ be a connected bipartite graph with bipartition $V (G)
= A \cup B$\,, $|A| = p$\,, $|B| = q$\,, $p+q=n$\,. Then $\rho (G)
\geq n - 2 + \sqrt{n^2 - 3pq}$ with equality if and only if $G\cong
K_{p,q}$\,. This was  shown by Das \cite{Da09} using a different
reasoning.  Note that $\rho (K_{p, q})=p + q - 2 + \sqrt{(p + q)^2 -
3pq}$  attains minimum if and only if $|p - q| \leq 1$\,. Therefore,
we have:

\begin{Theorem}
Among connected bipartite graphs with $n$ vertices, $K_{\lfloor n /
2 \rfloor, \lceil n / 2 \rceil}$ has minimum, while $P_n$ has
maximum distance spectral radius.
\end{Theorem}

Here we present a stronger lower bound for $\rho(G)$ for a bipartite
graph $G$\,, involving the maximum degrees in both partite sets.

\begin{Theorem} \label{biplow}
Let $G$ be a connected bipartite graph with
bipartition $V (G) = A \cup B$\,, $|A| = p$\,, $|B| = q$\,, $p+q=n$\,. Let
$\Delta_A$ and $\Delta_B$ be maximum degrees among vertices
from $A$ and $B$\,, respectively. Then
$$
\rho (G) \geq n - 2 + \sqrt{ n^2 - 4 pq + (3q - 2 \Delta_A) (3p - 2 \Delta_B) }
$$
with equality if and only if $G$ is a complete bipartite graph
$K_{p,q}$ or $G$ is a semi--regular graph with every vertex
eccentricity equal to $3$\,.
\end{Theorem}

\begin{Proof}
Let $A = \{1, 2, \ldots, p\}$ and $B = \{p + 1, p + 2, \ldots, p + q\}$\,. Let $x = (x_1, x_2, \ldots, x_n)^T$
be a Perron eigenvector of $D(G)$ corresponding to the largest eigenvalue $\rho (G)$\,, such that
$$
x_i = \min_{k \in A} x_k \quad \mbox{and} \quad x_j = \min_{k \in B} x_k\,.
$$
From the eigenvalue equation $\rho (G) \cdot x = D (G) \cdot x$\,,
written for the component $x_i$  we have
\begin{eqnarray*}
\rho(G) x_i &=& \sum_{k = 1}^p d_{ik} \cdot x_k + \sum_{k = p + 1}^{p + q} d_{ik} \cdot x_k \\[2mm]
&\geq & 2 (p - 1) x_i + d_i x_j + 3 (q - d_i) x_j \\[2mm]
&\geq& 2 (p - 1) x_i + (3 q - 2 \Delta_A) x_j\,.
\end{eqnarray*}
Analogously for the component $x_j$ we have
\begin{eqnarray*}
\rho(G) x_j &=& \sum_{k = 1}^p d_{jk} \cdot x_k + \sum_{k = p + 1}^{p + q} d_{jk} \cdot x_k \\[2mm]
&\geq & d_j x_i + 3 (p - d_j) x_i + 2 (q - 1) x_j \\[2mm]
&\geq& (3 p - 2 \Delta_B) x_i + 2 (q - 1) x_j\,.
\end{eqnarray*}
Combining these two inequalities, it follows
$$
(\rho (G) - 2 (p - 1))(\rho (G) - 2 (q - 1)) x_i x_j \geq (3q - 2 \Delta_A)(3p - 2 \Delta_B) x_i x_j\,.
$$
Since $x_k > 0$ for $1 \leq k \leq p + q$\,,
$$
\rho^2 (G) - 2 (p + q - 2)\rho (G) + 4 (p-1)(q-1) - (3q - 2
\Delta_A)(3p - 2 \Delta_B) \geq 0\,.
$$
From this inequality, we get the result.

For the case of equality, we have $x_i = x_k$ for $k = 1, 2, \ldots,
p$ and $x_j = x_k$ for $k = p + 1, p + 2, \ldots, p + q$\,. This
means that eigenvector $x$ has at most two different coordinates,
the degrees of vertices in $A$ are equal to $\Delta_A$\,, and the
degrees of vertices in $B$ are equal to $\Delta_B$\,, implying that
$G$ is a semi--regular graph. If $G$ is not a complete bipartite
graph, it follows from $p \Delta_A = q \Delta_B$ that $\Delta_A < q$
and $\Delta_B < p$ and the eccentricity of every vertex must be
equal to $3$\,.
\end{Proof}

It is evident that the lower bound in previous theorem improves the
bound in \cite{Da09} mentioned above.

Let $G$ be a connected bipartite graph with $n$ vertices,
diameter $d$ and bipartition $V (G) = A \cup B$\,, $|A| = p$\,, $|B| =
q$\,, $p+q=n$\,. Das in \cite{Da09} proved that
$$
\rho (G) \leq \frac{1}{2} \left [ d (n - 2) + \sqrt{d^2n^2 - 4pq (2d - 1)} \,\right ]
$$
for even $d$\,, and
\begin{eqnarray*}
\rho (G) &\leq & \frac{1}{2}  (d - 1) (n - 2)\\[2mm]
&&+ \ \frac{1}{2}\sqrt{(d - 1)^2n^2 +4\delta^2(d-1)^2- 4 pq(2d-1) - 4 d (d - 1)\delta n}
\end{eqnarray*}
for odd $d$\,.
Here we improve this result.

\begin{Theorem} \label{upev}
Let $G$ be a connected bipartite graph with $n$ vertices,
diameter $d$ and bipartition $V (G) = A \cup B$\,, $|A| = p$\,,
$|B| = q$\,, $p+q=n$\,.  Let $\delta_A$ and $\delta_B$ be the minimum degrees among vertices
from $A$ and $B$\,, respectively.
Then
\begin{eqnarray*}
\rho (G) & \leq & \frac{d}{2} \left ( n - 1 - \frac{d}{2} \right)\\[2mm]
&&+ \ \frac{1}{2} \sqrt{d^2 n^2 + 4 \delta_A \delta_B (d-2)^2 - 4 p q (2d - 1) - 4 (d-1)(d-2)(p\delta_A+q\delta_B)}
\end{eqnarray*}
for even $d$\,, and
\begin{eqnarray*}
\rho (G) &\leq& \frac{2 (d-1)n + 1-d^2}{4}\\[2mm]
&& + \ \frac{1}{2} \sqrt{(d - 1)^2 n^2 + 4 \delta_A \delta_B (d-1)^2 + 4 pq (2d - 1) - 4d (d - 1)(p\delta_A  + q\delta_B)}
\end{eqnarray*}
for odd $d$\,.
\end{Theorem}

\begin{Proof}
Let $A = \{1, 2, \ldots, p\}$ and $B = \{p + 1, p + 2, \ldots, p + q\}$\,. Let $x = (x_1, x_2, \ldots, x_n)^T$
be a Perron eigenvector of $D(G)$ corresponding to the largest eigenvalue $\rho (G)$\,, such that
$$
x_i = \max_{k \in A} x_k \quad \mbox{and} \quad x_j = \max_{k \in B} x_k\,.
$$
Suppose $d$ is even. From the eigenvalue equation $\rho (G) \cdot x = D (G) \cdot x$\,, written for the component $x_i$ we have
\begin{eqnarray*}
\rho(G) x_i &=& \sum_{k = 1}^p d_{ik} \cdot x_k + \sum_{k = p + 1}^{p + q} d_{ik} \cdot x_k \\[2mm]
&\leq& x_i \cdot \left[2 + 4 + \cdots + (d - 2)\right] + x_i \cdot d \left [p - 1 - \left (\frac{d}{2} - 1 \right ) \right ] \\[2mm]
&& + \ x_j \cdot \left[1 \cdot \delta_A + 3 + \cdots + (d - 3)\right] \\[2mm]
&& + \ x_j \cdot
(d - 1)\left [q - \left (\frac{d}{2} - 1 \right ) - (\delta_A -
1)\right ]\\[2mm]
&=& \left(-\frac{d^2}{4}-\frac{d}{4}+d p\right)
x_i+\left[(d-1)q-(d-2)\delta_A-\frac{d^2}{4}+\frac{3d}{2}-2\right]x_j\,,
\end{eqnarray*}
or equivalently
\begin{eqnarray*}
\left ( \rho(G) + \frac{d^2}{4} + \frac{d}{2} - d p \right ) x_i
&\leq& \left[(d-1)q-(d-2)\delta_A-\frac{d^2}{4}+\frac{3d}{2}-2\right]x_j\,.
\end{eqnarray*}
Note that for $d$ being an even number,
$-\frac{d^2}{4}+\frac{3d}{2}-2\le 0$\,.
We get
\begin{eqnarray*}
\left ( \rho(G) + \frac{d^2}{4} + \frac{d}{2} - d p \right ) x_i
&\leq& \left[ (d - 1)q  - (d - 2) \delta_A \right] x_j\,.
\end{eqnarray*}
Analogously for the component $x_j$ we have
\begin{eqnarray*}
\left ( \rho(G) + \frac{d^2}{4} + \frac{d}{2} - d q \right ) x_j
&\leq& \left[(d-1)p-(d-2)\delta_B\right]x_i\,.
\end{eqnarray*}
Combining these two inequalities, it follows
\begin{eqnarray*}
0 &\geq& \rho^2 (G) + \left [ \frac{d^2}{2} + d - d (p + q) \right] \rho (G) +\frac{d^2}{16} (2 + d - 4 p) (2 + d - 4 q) \\[2mm]
&&- \left[(d-1)q-(d-2)\delta_A\right]\left[(d-1)p-(d-2)\delta_B\right]\,.
\end{eqnarray*}
By analyzing the quadratic inequality and using $p + q = n$\,, we get
\begin{eqnarray*}
\rho (G) & \leq & \frac{d}{2} \left ( n - 1 - \frac{d}{2} \right) \\[2mm]
&&+ \ \frac{1}{2} \sqrt{d^2 n^2 + 4 \delta_A \delta_B (d-2)^2 - 4 p q (2d - 1) - 4 (d-1)(d-2)(p\delta_A+q\delta_B)}\,,
\end{eqnarray*}
as desired for even $d$\,.

Now suppose that $d$ is odd.
From the eigenvalue equation $\rho (G) \cdot x = D (G) \cdot x$\,,
written for component $x_i$ we have
\begin{eqnarray*}
\rho(G) x_i &=& \sum_{k = 1}^p d_{ik} x_k + \sum_{k = p + 1}^{p + q} d_{ik} x_k
\\[2mm]
&\leq& x_i \cdot [2 + 4 + \ldots + (d - 3)] + x_i \cdot (d-1) \left (p - 1 - \frac{d-3}{2}  \right )
\\[2mm]
&&+ x_j \cdot [1 \cdot \delta_A + 3 + \ldots + (d - 2)] + x_j \cdot
d\left (q - \frac{d-3}{2}  - \delta_A \right )\\[2mm]
&=& \left[-\frac{d^2}{4}+\frac{1}{4}+(d-1) p\right]
x_i+\left[dq-(d-1)\delta_A-\frac{d^2}{4}+d-\frac{3}{4}\right]x_j\,,
\end{eqnarray*}
or equivalently
\begin{eqnarray*}
\left [ \rho(G)+\frac{d^2}{4}-\frac{1}{4}-(d-1) p \right ] x_i
&\leq&
\left[dq-(d-1)\delta_A-\frac{d^2}{4}+d-\frac{3}{4}\right]x_j\,.
\end{eqnarray*}
Note that for $d$ being an odd number,
$-\frac{d^2}{4}+d-\frac{3}{4}\le 0$\,. We get
\begin{eqnarray*}
\left [\rho(G)+\frac{d^2}{4}-\frac{1}{4}-(d-1) p \right ] x_i &\leq&
\left [dq-(d-1)\delta_A\right ]x_j\,.
\end{eqnarray*}
Analogously for the component $x_j$ we have
\begin{eqnarray*}
\left [ \rho(G)+\frac{d^2}{4}-\frac{1}{4}-(d-1) q \right ] x_j
&\leq& \left [dp-(d-1)\delta_B\right ]x_i\,.
\end{eqnarray*}
Then the result for odd $d$ follows easily.
\end{Proof}


If the upper bound with even $d$ is attained in Theorem \ref{upev}
for $d=4$, then we have equal values of eigencomponents in both
partite sets, from which $G$ is semi--regular, for every vertex $v$
there is a unique vertex at distance $2$ and thus all vertices have
degree at most $2$, which is impossible. Thus the upper bound is
attained for even $d$ in Theorem \ref{upev} if and only if $d=2$,
$\delta_A=q, \delta_B=p$ and $G\cong K_{p,q}$\,. The upper bound is
attained for odd $d$ in Theorem \ref{upev} if and only if $d=1$,
$p=q=\Delta_A=\Delta_B=1$ and $G\cong K_{1,1}$, or $d=3$\,, $G$ is
semi--regular, any two vertices from the same partite set are at
distance $2$ and all vertex eccentricities are equal to $3$\,.

Under the conditions of Theorem \ref{upev}, let $\delta$ be the
minimum degree. Then $\delta=\min\{\delta_A, \delta_B\}$ and by
Theorem \ref{upev}, for even $d$ we have
\begin{eqnarray*}
\rho (G) &\leq & \frac{d}{2} \left ( n- 1 - \frac{d}{2} \right)\\[2mm]
&& + \frac{1}{2} \sqrt{  d^2 n^2 + 4 \delta^2 (d-2)^2 - 4 p q (2d - 1) - 4 \delta (d-1)(d-2) n}\,,
\end{eqnarray*}
 while for odd $n$ we have
\begin{eqnarray*}
\rho (G) &\leq& \frac{2 (d-1)n + 1-d^2}{4} \\[2mm]
&& + \frac{1}{2} \sqrt{(d - 1)^2 n^2 + 4 \delta^2 (d-1)^2 + 4 pq (2d
- 1) - 4\delta d (d - 1)n}\,.
\end{eqnarray*}
These two upper bounds for $\rho(G)$ improve the upper bounds for
$\rho(G)$ in \cite{Da09} mentioned above for even $d$ and for odd
$d$, respectively.

\section*{\normalsize 5.  LOWER BOUNDS FOR $DE$}

Let $G$ be a connected graph with $n\ge 2$ vertices. Note that $\rho (G)>0$\,.
Then
$$
DE(G)\geq 2\rho (G)
$$
with equality if and only if $G$ has exactly one positive
$D$-eigenvalue.  Thus, the lower bounds for $\rho (G)$ may  be converted
to lower bounds for $DE$\,.

Let $G$ be a connected graph with $n\ge 2$ vertices. Let $D_i$ be
the $i$-th row sum of ${\bf D}(G)$\,, i.~e.,
$D_i=\sum\limits_{j=1}^n d_{ij}$\,, where $i=1,2,\ldots,n$\,. It was
shown in \cite{ZhTr07} that
\begin{equation}
\label{bas} \rho (G) \geq \sqrt{\frac{1}{n}\,\sum_{i=1}^n D_i^2}
\end{equation}
with equality if and only if ${\bf D}(G)$ has equal row sums. (In
the view of matrix theory, by considering the powers of ${\bf D}(G)$\,, this
lower bound could be further improved.) Thus, we have:

\begin{Theorem} Let $G$ be a connected graph with $n\ge 2$
vertices. Then
$$
DE(G) \geq 2\,\sqrt{\frac{1}{n}\,\sum_{i=1}^n D_i^2}
$$
with equality if and only if $G$ has exactly one positive
$D$-eigenvalue and ${\bf D}(G)$ has equal row sums.
\end{Theorem}

Recall that trees \cite{CDGT}, connected unicyclic graphs \cite{BaKi05},
and $K_n$ have
exactly one positive $D$-eigenvalue. A complete
characterization of such graphs seems to be not known.

The Wiener index \cite{Ho,DoEn01} of a connected graph $G$ is defined as
$W(G)=\sum\limits_{i<j}d_{ij}=\frac{1}{2}\sum\limits_{i=1}^nD_i$\,.
From (\ref{bas}) and using the Cauchy--Schwarz inequality, as in
\cite{ZhTr07}, we get
$$
\rho (G) \geq \frac{2W(G)}{n}
$$
with equality if and only if ${\bf D}(G)$ has equal row sums.
Thus, for $m$ being the number of edges of $G$\,,
\[
\rho (G)\ge 2(n-1)-\frac{2m}{n}
\]
with equality if and only if $G \cong K_n$ or $G$ is a regular graph of
diameter two. It follows

\begin{Theorem} \label{P2}
Let $G$ be a connected graph with $n\geq 2$
vertices and $m$ edges. Then
\[
DE(G)\ge \frac{4W(G)}{n}
\]
with equality if and only if ${\bf D}(G)$ has equal row sums and $G$ has exactly one positive $D$-eigenvalue.
Moreover,
\[
DE(G)\ge 4(n-1)-\frac{4m}{n}
\]
with equality if and only if either  $G \cong K_n$ or  $G$ is a regular
graph of diameter two with exactly one positive $D$-eigenvalue.
\end{Theorem}

Ramane et al.~\cite{RRGRAW} conjectured that among the $n$-vertex
connected graphs, the complete graph $K_n$ is the unique graph with
the smallest distance energy (equal to $2(n-1)$). For a connected
graph $G$ with $n$ vertices and $m$ edges, $2m\le n(n-1)$ with
equality if and only if $G \cong K_n$\,. By Theorem \ref{P2},
this conjecture is true. A direct reasoning is as follows: Note that
$K_n$ for $n\ge 2$ has exactly one positive $D$-eigenvalue. From
(\ref{eq-addedge}), we have $\rho_1\ge n-1$\,, and then $DE(G)\ge
2\rho_1\ge 2(n-1)$ with equalities if and only if $G \cong K_n$\,.


Let $G$ be a graph. The line graph $L(G)$ of $G$ has the edges of
$G$ as vertices, and vertices of $L(G)$ are adjacent if the
corresponding edges of $G$ have a vertex in common. The
cocktail party graph $CP(a)$ is the graph obtained by deleting $a$
disjoint edges from the complete graph $K_{2a}$\,. Thus, $CP(a)$ is
a regular graph of degree $2a-2$\,.

Let $G$ be a graph with vertex set $\{v_1, v_2, \dots, v_n\}$\,,
and let $a_1,a_2,\ldots,a_n$ be nonnegative integers. The generalized line graph
$L(G;a_1,a_2 \ldots,a_n)$ consists of the disjoint union of the line
graph $L(G)$ and the cocktail party graphs $CP(a_1),CP(a_2),\ldots,CP(a_n)$\,,
together with all edges joining a vertex $\{v_i,v_j\}$ of
$L(G)$ with each vertex of $CP(a_i)$ and $CP(a_j)$\,.

A regular graph $G$ of diameter two has exactly one positive
$D$-eigenvalue if and only if $\lambda_n\geq -2$\,. For a
generalized line graph, its least eigenvalue is at least $-2$\,.  An
exceptional graph is a connected graph, other than a generalized
line graph, with least eigenvalue at least $-2$\,. From \cite{CRS},
a graph $G$ is a regular graph of diameter two with $\lambda_n\ge
-2$ if and only if $G$ is a cocktail party graph, or $G$ is a
regular line graph of diameter two (equivalently, $G$ is the line
graph with diameter two of a regular graph or of a semi--regular
graph), or $G$ is a regular exceptional graph of diameter two. A
list of all $187$ regular exceptional graphs is given in Table A3 of
\cite{CRS}. These graphs are listed in such a way that it is not
possible to see what their diameters are. Of these, exactly $7$
graphs are strongly regular, and thus have diameter two. However, it
may be that there are other, not strongly regular graphs with
diameter two. We do not attempt to count them.

\begin{Corollary}
Let $G$ be a connected graph with $n\ge 2$ vertices and $m$ edges.
Then
$$
DE(G) \geq 4(n-1)-\frac{4m}{n}
$$
with equality if and only if  $G \cong K_n$\,, or  $G$ is the cocktail party graph,
or $G$ is a regular line graph of diameter two, or $G$ is a regular exceptional
graph of diameter two.
\end{Corollary}

By discussion above and Theorem \ref{low1}, we have

\begin{Corollary}
Let $G$ be a connected graph with $n$ vertices, maximum degree $\Delta_1$ and second maximum degree $\Delta_2$\,.
Then
$$
DE (G) \geq 2 \sqrt{(2n-2-\Delta_1)(2n-2-\Delta_2)}
$$
with equality if and only if $G \cong K_n$\,, or  $G$ is the cocktail party graph,
or $G$ is a regular line graph of diameter two, or $G$ is a regular exceptional
graph of diameter two.
\end{Corollary}

For the graph $G$\,, the first Zagreb index of $G$ is defined as
$M_1(G)=\sum\limits_{i=1}^nd_i^2$\, [29--32]. 

Let $G$ be a triangle-- and quadrangle--free connected graph with $n\ge
2$ vertices and $m$ edges.  Then from \cite{ZhTr07} it follows
\[
\rho (G)\ge 3(n-1)-\frac{2m}{n}-\frac{M_1(G)}{n}
\]
with equality if and only if ${\bf D}(G)$ has equal  row sums and
the diameter of $G$  is at most three, and thus
\[
DE(G)\ge 2\left[3(n-1)-\frac{2m}{n}-\frac{M_1(G)}{n}\right]
\]
with equality if and only if $G$ has exactly one positive
$D$-eigenvalue, ${\bf D}(G)$ has equal row sums and the diameter of
$G$  is at most three.

Let $G$ be a connected bipartite graph with $p$ vertices in one
partite set and $q$ vertices in the other partite set. Recall that
$\rho (G)\ge p+q-2+\sqrt{p^2+q^2-pq}$ with equality if and only if
$G\cong K_{p,q}$ and that $K_{p,q}$ has exactly one positive
$D$-eigenvalue if and only if $3pq\le 4(n-1)$. This implies:

\begin{Theorem} Let $G$ be a connected bipartite graph with $p$ vertices in one
partite set and $q$ vertices in the other partite set. Then
\[
DE(G)\ge 2\left(p+q-2+\sqrt{p^2+q^2-pq}\,\right)
\]
with equality if and only if $G \cong K_{p,q}$ with  $3pq\le 4(n-1)$\,.
\end{Theorem}

From this theorem, we have: if $G$ is a connected bipartite
graph with $n\ge 2$ vertices, then
\[
DE(G)\ge 2\left(n-2+\sqrt{n^2-3\left\lfloor\frac{n}{2}\right\rfloor
\left\lceil\frac{n}{2}\right\rceil}\,\right)
\]
with equality if and only if $G \cong
K_{\lfloor n/2\rfloor, \lceil n/2\rceil}$ with  $n=2,3,4$\,.

Let $G$ be a connected bipartite graph with
bipartition $V (G) = A \cup B$\,, $|A| = p$\,, $|B| = q$\,, $p+q=n$\,. Let
$\Delta_A$ and $\Delta_B$ be maximum degrees among vertices
from $A$ and $B$\,, respectively.  By Theorem \ref{biplow} and the discussion above, we have:
$$
DE(G) \geq 2(n - 2) + 2\sqrt{ n^2 - 4 pq + (3q - 2 \Delta_A) (3p - 2 \Delta_B) }
$$
with equality if and only if $G$ is a complete bipartite graph
$K_{p,q}$ or $G$ is a semi--regular graph with every vertex
eccentricity equal to $3$ and exactly one positive $D$-eigenvalue.

\begin{Theorem}
Let $G$ be a connected graph with $n$ vertices. If $\overline{G}$ is
also connected, then
$$
DE(G)+DE(\overline{G}) \geq 6(n-1)
$$
with equality if and only if $G$ and $\overline{G}$ both have
exactly one positive $D$-eigenvalue and are both regular graphs of
diameter two.
\end{Theorem}

\begin{Proof}
By Theorem \ref{P2},
\[
DE(G)+DE(\overline{G})\ge 8(n-1)-\frac{2n(n-1)}{n}=6(n-1)
\]
with equality if and only if $G$ and $\overline{G}$ both have
exactly one positive $D$-eigenvalue and are both regular graphs of
diameter two.
\end{Proof}

\section*{\normalsize 6.  UPPER BOUNDS FOR $DE$}

In the following we discuss upper bounds for the distance energy of graphs of diameter at most two.

Let $G$ be a simple graph with $n$ vertices. The energy of the graph
$G$ is defined as \cite{Gut1,Gut2}
\[
E(G)=\sum_{i=1}^n |\lambda_i|\,.
\]

The singular eigenvalues of a (complex) matrix ${\bf X}$ are the square
roots of the eigenvalues of the matrix ${\bf XX}^*$\,, where ${\bf
X}^*$ denotes the conjugate transpose of the matrix ${\bf X}$\,.
For an $n\times n$ matrix ${\bf X}$\,, its singular values are denoted by
$s_i({\bf X})$\,,  $i=1, 2, \dots, n$\,. Then
$$
E(G)=\sum_{i=1}^ns_i({\bf A}(G))
$$
$$
DE(G)=\sum_{i=1}^ns_i({\bf D}(G))\,.
$$

\begin{Lemma} \textnormal{\cite{Fan}}
\label{Fan} Let ${\bf X}$ and ${\bf Y}$ be  $n\times n$ matrices. Then
 $\sum\limits_{i=1}^ns_i({\bf X}+{\bf Y})\le \sum\limits_{i=1}^ns_i({\bf
X})+\sum\limits_{i=1}^ns_i({\bf Y})$\,.
\end{Lemma}

\begin{Theorem}
Let $G$ be a connected graph with $n$ vertices and diameter at most two. Then
$$
DE(G)\le 2(n-1)+E(\overline{G})\,.
$$
\end{Theorem}

\begin{Proof} Note that
$$
{\bf D}(G)={\bf J}_n-{\bf I}_n+{\bf A}(\overline{G})\,.
$$
Let ${\bf X}={\bf J}_n-{\bf I}_n$ and ${\bf Y}={\bf A}(\overline{G})$ in Lemma \ref{Fan}, we have
$$
DE(G)\le 2(n-1)+E(\overline{G})\,,
$$
as desired.
\end{Proof}

In \cite{Ko}, it was shown that for a graph with $n$ vertices, its energy is bounded from above by $\frac{n}{2}(\sqrt{n}+1)$\,.
Thus, for the  connected graph $G$  with $n$ vertices and diameter at most two,
$$
DE(G)\le \frac{n}{2}(\sqrt{n}+1)+2(n-1)\,.
$$
For $n\ge 26$\,, this is better than the bound given in \cite{IGV}:
$$
DE(G)\leq \sqrt{2n(2n^{2}-2n-3m)}\,,
$$
where $m$ is the number of edges of $G$\,,
because for $n\ge 26$\,,
$$
\frac{n}{2}(\sqrt{n}+1)+2(n-1)<n\sqrt{n-1}\le \sqrt{2n(2n^{2}-2n-3m)}\,.
$$

In \cite{IGV}, it was shown that for a graph $G$ with $n$ vertices, $m$ edges and diameter two,
$$
DE(G)\le \frac{1}{n}(2n^2-2n-2m)+\frac{1}{n}\sqrt{(n-1)[(2n+m)(2n^2-4m)-4n^2]}\,.
$$
This upper bound for  $K_{1,n-1}$  is equal to
$2n-4+\frac{2}{n}+\frac{1}{n}\sqrt{6n^4-24n^3+34n^2-20n+4}$\,, while
the bound in the previous theorem is $DE(K_{1, n-1})\le
2(n-1)+2(n-2)=4n-6$\,. The latter is better than the former for
$n\ge 5$\,.

\bigskip

\noindent {\bf Added after publication}: Corollary 1 confirms Conjecture 3 in [G. Caporossi, E. Chasset, B. Furtula, Some conjectures and properties on distance energy, at: {\tt http://www.\\ gerad.ca/fichiers/cahiers/G-2009-64.pdf}].

\vspace{4mm}

\baselineskip=0.27in

\noindent {\it Acknowledgement.\/} This work was supported by the
Guangdong Provincial Natural Science Foundation of China (Grant
No.~8151063101000026) and by Research Grant 144007 of the Serbian
Ministry of Science and Technological Development. We thank
Professor Ivan Gutman for his help and encouragement.

\renewcommand{\refname}{
\end{document}